\documentclass[a4paper,12pt]{article}
\usepackage{amssymb}
\usepackage{amsthm}
\usepackage{tikz}
\usepackage{pgfplots}
\usetikzlibrary{decorations.markings}
\usepackage{amsmath}
\usepackage{todonotes}
\usepackage{url}
\usepackage[margin=2.5cm]{geometry}
\usepackage{color}
\usepackage{graphicx}
\usepackage{placeins,caption} 
\usepackage{multicol}   
\usepackage{comment}
\usepackage{hyperref}

\hypersetup{colorlinks=true, linkcolor=blue, citecolor=blue, urlcolor=blue}

\newcommand{\tr}[1]{\mathrm{Tr}\,#1}

\def\A{\mbox{\boldmath $A$}}
\def\B{\mbox{\boldmath $B$}}

\def\I{\mbox{\boldmath $I$}}

\let\lcm\relax 
\DeclareRobustCommand{\lcm}{\mathop{\operator@font lcm}}

\newtheorem{theorem}{Theorem}[section]
\newtheorem{lemma}[theorem]{Lemma}

\newtheorem{proposition}[theorem]{Proposition}

\newtheorem{conjecture}[theorem]{Conjecture}
\newtheorem{problem}[theorem]{Problem}

\theoremstyle{definition}

\pgfplotsset{compat=1.18}

\begin{document}
	
	\title{There are no excess one digraphs}
	
	\author{Slobodan Filipovski $^a$ \\
		\texttt{\footnotesize slobodan.filipovski@famnit.upr.si} \\
		{\footnotesize ORCID: 0000-0002-7286-4954}  \\
		\and 
		Arnau Messegu\'e $^{b}$ \\ \texttt{\footnotesize arnau.messegue@udl.cat}  \\
		{\footnotesize ORCID: 0000-0002-7425-7592}
		\and
		Josep M. Miret $^{b}$ \\ \texttt{\footnotesize josepmaria.miret@udl.cat} \\
		{\footnotesize ORCID: 0000-0003-4631-4294}
		\and 
		James Tuite $^{c,d}$ \\ \texttt{\footnotesize james.t.tuite@open.ac.uk} \\ 
		{\footnotesize ORCID: 
			0000-0003-2604-7491}
	}
	
	\maketitle
	
	\begin{center}
		$^a$ FAMNIT, University of Primorska, Koper, Slovenia \\
		
		$^b$ Departament de Matem\`{a}tica, Universitat de Lleida, Spain\\

		$^c$ School of Mathematics and Statistics, Open University, Milton Keynes, UK\\
		
		$^d$ Department of Informatics and Statistics, Klaip\.{e}da University, Lithuania\\
		
	\end{center}
	
	\begin{abstract}
		A digraph $G$ is \emph{$k$-geodetic} if for any pair $u,v \in V(G)$ there is at most one $u,v$-walk of length not exceeding $k$. The order of a $k$-geodetic digraph with minimum out-degree $d$ is bounded below by the directed Moore bound $M(d,k) = 1 + d + d^2+ \cdots +d^k$. It is known that the Moore bound cannot be achieved for $d,k \geq 2$. A $k$-geodetic digraph with minimum degree $d$ and order one greater than the Moore bound has \emph{excess one}. In this paper we prove a conjecture that no excess one digraphs exist for $d,k \geq 2$, thus complementing the result of Bannai and Ito on the non-existence of undirected graphs with excess one.
	\end{abstract}
	
	\noindent
	{\bf Keywords:} digraph, geodecity, excess, degree/diameter, cage
	
	\noindent
	AMS Subj.\ Class.\ (2020): 05C12, 05C20, 05C35
	\section{Introduction}
	
	A digraph $G = (V,A)$ consists of a set $V(G)$ of vertices and a set $A(G)$ of directed arcs, i.e.\ ordered pairs $(u,v)$ of vertices. We write $u \rightarrow v$ to indicate the existence of an arc $(u,v)$. In a digraph we do not allow loops $u \rightarrow u$ or multiple arcs $u \rightarrow v$ with the same direction; if such connections are present, the network is a \emph{pseudodigraph}. The \emph{out-neighbourhood} $N^+(u)$ of a vertex $u$ is $\{ v \in V(G):u \rightarrow v\} $ and the \emph{out-degree} is $d^+(u) = |N^+(u)|$. Similarly, the \emph{in-neighbourhood} $N^-(u)$ of $u$ is $\{ v \in V(G):v \rightarrow u\} $ and the \emph{in-degree} is $d^-(u) = |N^-(u)|$. The digraph is \emph{out-regular} if $d^+(u) = d$ for some $d$ and all $u \in V(G)$ and is \emph{diregular} if $d^-(u) = d^+(u) = d$ for some $d$ and all $u \in V(G)$. A \emph{directed walk} of length $\ell $ is a sequence $u_0,u_1,\dots ,u_{\ell }$ of not necessarily distinct vertices, such that $u_i \rightarrow u_{i+1}$ for $0 \leq i \leq \ell -1$. The \emph{distance} $d(u,v)$ between two vertices $u,v \in V(G)$ is the length of a shortest directed walk from $u$ to $v$, when such a walk exists. Note that we can have $d(uv) \neq d(v,u)$. The diameter of $G$ is $\max \{ d(u,v):u,v \in V(G)\} $. A digraph is \emph{strongly connected} if for any pair $u,v$ there is a directed walk in $G$ from $u$ to $v$.

\begin{figure}[h]
		\centering
		\begin{tikzpicture}[x=0.4mm,y=-0.4mm,inner sep=0.2mm,scale=0.25,very thick,vertex/.style={circle,draw,minimum size=10,fill=lightgray}]
			\node at (181.4,173.2) [vertex] (v1) {};
			\node at (222.6,300.4) [vertex] (v2) {};
			\node at (317.2,211.5) [vertex] (v3) {};
			\node at (95.9,379.5) [vertex] (v4) {};
			\node at (679.4,379.5) [vertex] (v5) {};
			\node at (387.6,87.8) [vertex] (v6) {};
			\node at (387.6,671.2) [vertex] (v7) {};
			\node at (118.1,267.9) [vertex] (v8) {};
			\node at (118.1,491.1) [vertex] (v9) {};
			\node at (657.2,267.9) [vertex] (v10) {};
			\node at (657.2,491.1) [vertex] (v11) {};
			\node at (499.3,110) [vertex] (v12) {};
			\node at (276,110) [vertex] (v13) {};
			\node at (499.3,649) [vertex] (v14) {};
			\node at (276,649) [vertex] (v15) {};
			\node at (181.4,585.8) [vertex] (v16) {};
			\node at (593.9,173.2) [vertex] (v17) {};
			\node at (593.9,585.8) [vertex] (v18) {};
			\node at (547.6,297.7) [vertex] (v19) {};
			\node at (456.8,212.6) [vertex] (v20) {};
			\path
			(v1) edge[->] (v2)
			(v1) edge[->] (v3)
			(v8) edge[->] (v1)
			(v13) edge[->] (v1)
			(v2) edge[->] (v4)
			(v2) edge[->] (v5)
			(v18) edge[->] (v2)
			(v3) edge[->] (v6)
			(v3) edge[->] (v7)
			(v18) edge[->] (v3)
			(v4) edge[->] (v8)
			(v4) edge[->] (v9)
			(v19) edge[->] (v4)
			(v5) edge[->] (v10)
			(v5) edge[->] (v11)
			(v19) edge[->] (v5)
			(v6) edge[->] (v12)
			(v6) edge[->] (v13)
			(v20) edge[->] (v6)
			(v7) edge[->] (v14)
			(v7) edge[->] (v15)
			(v20) edge[->] (v7)
			(v8) edge[->] (v12)
			(v15) edge[->] (v8)
			(v13) edge[->] (v9)
			(v9) edge[->] (v14)
			(v9) edge[->] (v16)
			(v10) edge[->] (v13)
			(v14) edge[->] (v10)
			(v10) edge[->] (v17)
			(v12) edge[->] (v11)
			(v11) edge[->] (v15)
			(v11) edge[->] (v18)
			(v12) edge[->] (v17)
			(v14) edge[->] (v18)
			(v15) edge[->] (v16)
			(v16) edge[->] (v19)
			(v16) edge[->] (v20)
			(v17) edge[->] (v19)
			(v17) edge[->] (v20)
			;
		\end{tikzpicture}
		\quad
		\begin{tikzpicture}[x=0.4mm,y=-0.4mm,inner sep=0.2mm,scale=0.3,very thick,vertex/.style={circle,draw,minimum size=10,fill=lightgray}]
			\node at (532,328) [vertex] (v1) {};
			\node at (403,75) [vertex] (v2) {};
			\node at (56,252) [vertex] (v3) {};
			\node at (256,52) [vertex] (v4) {};
			\node at (79,399) [vertex] (v5) {};
			\node at (509,181) [vertex] (v6) {};
			\node at (332,528) [vertex] (v7) {};
			\node at (185,505) [vertex] (v8) {};
			\node at (464,460) [vertex] (v9) {};
			\node at (124,120) [vertex] (v10) {};
			\node at (332,52) [vertex] (v11) {};
			\node at (532,252) [vertex] (v12) {};
			\node at (256,528) [vertex] (v13) {};
			\node at (509,399) [vertex] (v14) {};
			\node at (79,181) [vertex] (v15) {};
			\node at (56,328) [vertex] (v16) {};
			\node at (464,120) [vertex] (v17) {};
			\node at (403,505) [vertex] (v18) {};
			\node at (185,75) [vertex] (v19) {};
			\node at (123,460) [vertex] (v20) {};
			\path
			(v1) edge[->] (v11)
			(v1) edge[->] (v12)
			(v1) edge[<-] (v13)
			(v1) edge[<-] (v14)
			(v2) edge[->] (v11)
			(v2) edge[->] (v12)
			(v2) edge[<-] (v17)
			(v2) edge[<-] (v18)
			(v3) edge[<-] (v11)
			(v3) edge[->] (v13)
			(v3) edge[<-] (v15)
			(v3) edge[->] (v16)
			(v4) edge[<-] (v11)
			(v4) edge[<-] (v15)
			(v4) edge[->] (v18)
			(v4) edge[->] (v19)
			(v5) edge[<-] (v12)
			(v5) edge[->] (v14)
			(v5) edge[<-] (v16)
			(v5) edge[->] (v20)
			(v6) edge[<-] (v12)
			(v6) edge[->] (v15)
			(v6) edge[<-] (v16)
			(v6) edge[->] (v17)
			(v7) edge[<-] (v13)
			(v7) edge[<-] (v14)
			(v7) edge[->] (v18)
			(v7) edge[->] (v19)
			(v8) edge[->] (v13)
			(v8) edge[->] (v16)
			(v8) edge[<-] (v19)
			(v8) edge[<-] (v20)
			(v9) edge[->] (v14)
			(v9) edge[<-] (v17)
			(v9) edge[<-] (v18)
			(v9) edge[->] (v20)
			(v10) edge[->] (v15)
			(v10) edge[->] (v17)
			(v10) edge[<-] (v19)
			(v10) edge[<-] (v20)
			;
		\end{tikzpicture}
		\caption{Geodetic cages for $d=2,k=3$ with excess $\epsilon = 5$}
		\label{fig:d2k3n20}
	\end{figure}
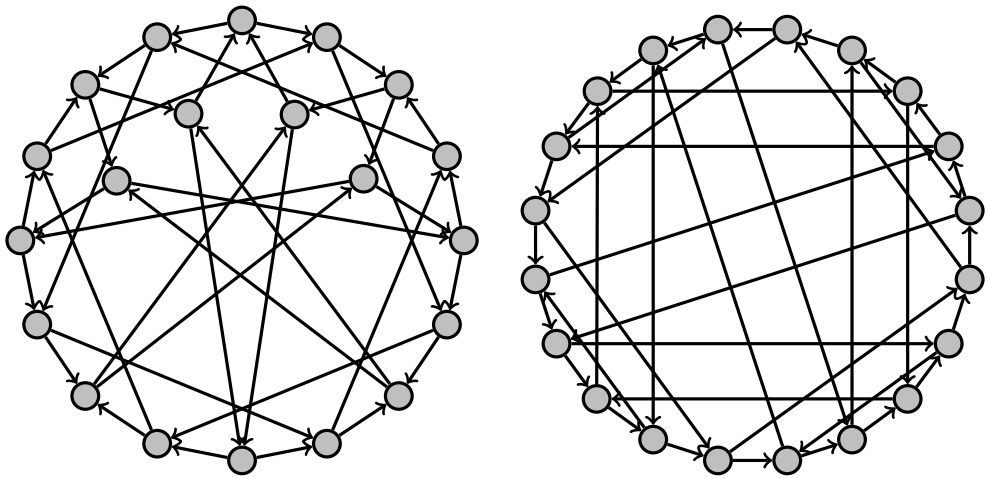

	The classical degree/girth problem asks for the smallest $d$-regular graphs with given girth $g$ (see~\cite{survey} for a survey). This problem was generalised to the setting of directed graphs in~\cite{Sil} as follows. A digraph $G$ is \emph{$k$-geodetic} if for any $u,v \in V(G)$ there is at most one directed $u,v$-walk in $G$ with length not exceeding $k$. The \emph{degree/geodecity problem} consists of finding the smallest possible order $N(d,k)$ of a $k$-geodetic digraph with minimum out degree $d$. The extremal graphs are known as \emph{geodetic cages}. Two geodetic cages are shown in Figure~\ref{fig:d2k3n20}. 
	
	The $k$-geodecity property implies that, for any vertex $u$ of a $k$-geodetic digraph, all vertices in the breadth-first search tree, or \emph{Moore tree}, of depth $k$ rooted at $u$ are distinct. It follows that the order of such a digraph is bounded below by the \emph{directed Moore bound} \[ N(d,k) \geq M(d,k) = 1+d+d^2+\cdots +d^k.\]  The existence of geodetic cages for all $d,k$ was proven in~\cite{networks}.
	
	If a digraph $G$ meets the Moore bound $M(d,k)$, then $G$ must have diameter $k$ and would be a \emph{directed Moore digraph}. It is shown in~\cite{BriTou,PZ74} that no directed Moore graphs exist, except for the trivial case of complete digraphs and directed $(k+1)$-cycles. Hence in general a geodetic cage $G$ will have order $N(d,k) = M(d,k) + \epsilon$, where $\epsilon $ is the \emph{excess} of $G$. A $k$-geodetic digraph with minimum out-degree $d$ and excess $\epsilon $ will be referred to as a \emph{$(d,k;+\epsilon )$-digraph}.
	
	A natural question is whether there exist digraphs with excess one, apart from the uninteresting case of directed $(k+2)$-cycles. In such a digraph, for any vertex $u$ there is a unique vertex $o(u)$, called the \emph{outlier} of $u$, that has distance $> k$ from $u$. It is shown in~\cite{MilMirSil} that the function $o$ is a digraph automorphism, which, by construction, is fixed-point free.
	
	It is known that any excess one digraph must be diregular~\cite{MilMirSil}. It was shown in~\cite{Sil} that there are no digraphs with excess one for degree $d = 2$ and in the two articles~\cite{MilMirSil,Tuite2024} non-existence was also shown for $k \leq 4$ when $d > 1$. It was conjectured in~\cite{thesis} that there are no excess one digraphs for $d,k \geq 2$. Digraphs with larger excess are discussed in the articles~\cite{d2e2nondiregular,regularity,d2e3,d2e2diregular,networks}. We mention also that the degree/geodecity problem has many similarities with the associated \emph{directed degree/diameter problem}; in particular, the analogue of excess one digraphs with an outlier automorphism $o$ in the degree/diameter problem is digraphs with defect one, with repeat automorphism $r:V(G) \rightarrow V(G)$; we refer the reader to more details of the problem in~\cite{beyond}.

	In this article we prove the conjecture of~\cite{thesis} by showing that there are no digraphs with excess one for $d,k \geq 2$. This complements the result of Bannai and Ito~\cite{BannaiIto} that there are no undirected graphs with excess one (apart from cycles).

	The plan of the paper is as follows. In Section~\ref{sec:pseudodigraph} we associate with an excess one digraph $G$ a quotient pseudodigraph $\tilde{G}$ with vertex set consisting of the orbits of the outlier function on $G$, and analyse the spectrum of $\tilde{G}$ in order to bound the number of orbits. We thereby prove the non-existence of excess one digraphs for all save a finite number of combinations of $d$ and $k$. In Section~\ref{sec:nonexistence} we rule out the existence of excess one digraphs for the exceptional combinations of $d$ and $k$. Finally, in Section~\ref{sec:conclusion} we reflect on the implications of this result, and suggest some directions for future research.

	\section{Non-existence for large $d$ and $k$}\label{sec:pseudodigraph}
	
	Let $G$ represent a given $(d,k;+1)$-digraph with degree $d \geq 2$. We will denote the order of $G$ by $N = M(d,k)+1$. In this section, we introduce a quotient pseudodigraph $\tilde{G}$ derived from $G$, which we call the \emph{orbital}-\emph{connection digraph}. Recall that the outlier function $o$ of $G$ is a digraph automorphism. We now consider the orbits of the outlier automorphism on $G$. For $v \in V(G)$, the \emph{outlier orbit} of $v$ is
	$$
	\mathcal{O}(v)=\left\{ v,o(v), \ldots, o^{ord_o(v)-1}(v)\right\} ,
	$$
	where $ord_o(v)$ is the \emph{order} of $v$ under the action of $o$, i.e.\ $ord_o(v)$ is the smallest positive integer $t$ such that $o^t(v) = v$. We denote the smallest order of a vertex under the outlier automorphism in $G$ by $\omega (G)$ and call this the \emph{index} of $G$. Observe that every vertex of $G$ has order at least two (and so $\omega (G) \geq 2$), and it is shown in~\cite{Tuite2024} that some vertex must have order at least three. The \emph{permutation structure} of $G$ is defined to be \[ (m_1,m_2,\ldots,m_N),\] where, for $1 \leq j \leq N$, $m_j$ is the number of outlier orbits having cardinality $j$. Note that these numbers satisfy $N=\sum_{j=1}^N jm_j$ and $m_1 = 0$.

	
	To construct $\tilde{G}$, we define the vertex set of $\tilde{G}$ to be the set of outlier orbits of $G$, that is,
	$$
	V(\tilde{G})=\left\{\mathcal{O}_1,\mathcal{O}_2,\ldots,\mathcal{O}_w\right\}, \quad w=\sum_{j=1}^N m_j.
	$$
	For each orbit $\mathcal{O}_i$, we choose a vertex $v(\mathcal{O}_i)$ in $\mathcal{O}_i$ as a \emph{representative} of $\mathcal{O}_i$.
	Then the number of arcs in $\tilde{G}$ from $\mathcal{O}_i$ to $\mathcal{O}_j$ is defined to be the number of arcs from the vertex $v(\mathcal{O}_i)$ to the orbit $v(\mathcal{O}_j)$ in the base digraph $G$. In particular, if $i = j$, there will $t$ loops at the vertex $\mathcal{O}_i$ of $\tilde{G}$ if there are $t$ arcs in $G$ from $v(\mathcal{O}_i)$ to other vertices of the orbit $\mathcal{O}_i$. As each vertex $v(\mathcal{O}_i)$ has out-degree $d$ in $G$, the pseudodigraph $\tilde{G}$ is out-regular with degree $d$ (counting multiple arcs and loops). The pseudodigraph $\tilde{G}$ is the quotient digraph of $G$ over the equitable partition of $G$ given by the outlier orbits, see Chapter 9 of~\cite{GR}. It is easily seen that, up to isomorphism, the construction does not depend on the choice of representatives from the orbits. The following lemma follows directly from the $k$-geodetic property and the definition of the outlier function. 
	
	\begin{lemma}
		\label{properties}
		For all $j,h$ in the range $1\leq j,h \leq w$, it holds that: 
		\begin{itemize}
			\item[i)] if $h \neq j$, then, for each vertex $x$ in the orbit $\mathcal{O}_h$, $G$ contains exactly one directed walk with length at most $k$ from $v(\mathcal{O}_j)$ to $x$, and
			\item[ii)]  if $h = j$, then, for each vertex $y$ in $\mathcal{O}_j - \{ o(v(\mathcal{O}_j))\} $, $G$ contains exactly one directed walk with length at most $k$ from $v(\mathcal{O}_j)$ to $y$. If $y = v(\mathcal{O}_j)$, this walk has length zero, whilst no such walk exists if $y = o(v(\mathcal{O}_j))$.
		\end{itemize}
	\end{lemma}
	
	Using the properties in Lemma~\ref{properties} we can derive a matrix equation that must be satisfied by the adjacency matrix $\tilde{\A}$ of $\tilde{G}$.  
	
	\begin{proposition}\label{pro:main eqn}
		If $\tilde{\A}$ is the adjacency matrix of $\tilde{G}$, then:
		\begin{equation}
			\label{eqn:B}
			\I_w+\tilde{\A} + \tilde{\A}^2+\cdots +\tilde{\A}^k = \B -\I_w,
		\end{equation}
		where $\B$ is the matrix $(b_{ij})_{1\leq i,j\leq w}$ with entries defined by $b_{ij} = |\mathcal{O}_j|$, i.e.\ each entry in the $j$-th column of $B$ is the number of vertices in the orbit $\mathcal{O}_j$ of $G$.
	\end{proposition} 
	\begin{proof}
		The $(j,h)$-coefficient of the matrix
		$$\left(\I_w+\tilde{\A} + \tilde{\A}^2+\cdots+\tilde{\A}^k\right)$$
		counts the number of walks of length at most $k$ in $\tilde{G}$ from vertex $\mathcal{O}_j$ to vertex $\mathcal{O}_h$. This quantity equals the total number of walks of length at most $k$ in $G$ from the vertex $v(\mathcal{O}_j)$ to the vertices in the orbit $\mathcal{O}_h$. If $j \neq h$, then by part {\em i)} of Lemma~\ref{properties}, this number equals $|\mathcal{O}_h|$. If $j = h$, then by part {\em ii)} of Lemma~\ref{properties}, the number is $|\mathcal{O}_h|-1$. Therefore the above polynomial in $\tilde{\A }$ is equal to $\B -\I_w$. 
	\end{proof}
	We proceed to use Proposition~\ref{pro:main eqn} to derive information on the spectrum of $G$. The characteristic polynomial of the matrix $\B-\I_w$ of order $w$ is given in the following lemma.
	
	\begin{lemma}\label{lem:phi B-I}
		\begin{equation}
			\label{eqn:2}
			\phi(\B-\I_w,x) = (x-(N-1))(x+1)^{w-1}.
		\end{equation}
	\end{lemma}
	\begin{proof}
		The eigenvalues of the matrix $\B$ are $N$ with multiplicity one and $0$ with multiplicity $w-1$. Hence, the characteristic polynomial of $\B$ is
		$$
		\phi(\B,x) = (x-N)x^{w-1}.
		$$
		Since $\I_w+ (\B-\I_w)=\B$, the claim follows on substituting $x$ by $1+x$ in $\phi(\B,x)$.
	\end{proof}
	Lemma~\ref{lem:phi B-I} and Equation~\ref{eqn:B} yield an exact expression for the characteristic polynomial of $\tilde{\A}$.
	
	\begin{proposition}\label{prop:k divides w-1}
		If a $(d,k;+1)$-digraph $G$ exists, then $k \mid w-1$ and the characteristic polynomial of $\tilde{\A}$ is: 
		$$
		\phi (\tilde{\A},x)=(x-d)(2+x+\cdots +x^k)^{(w-1)/k}$$
	\end{proposition}
	
	\begin{proof}
		By (\ref{eqn:B}) and Lemma~\ref{lem:phi B-I}, $\tilde{\A}$ has an eigenvalue $d$ and the remaining eigenvalues satisfy $1+\lambda + \dots + \lambda ^k = -1$. As the polynomial $2+x+\cdots + x^k$ is irreducible over $\mathbb Q$, it follows that the characteristic polynomial of $\tilde{\A}$ is as stated, which implies that $k \mid (w-1)$.
	\end{proof}
	
	We now proceed to bound the number $w$ of orbits in $G$. By Proposition~\ref{prop:k divides w-1}, we can write $w = 1+ak$, where $a$ is a non-negative number. Proposition~\ref{prop:k divides w-1} implies the following identity for the trace of the powers of the adjacency matrix of $\tilde{G}$. 
	
	\begin{lemma}\label{lem:trace A l}
		For all $\ell$ with $1\leq \ell\leq k-1$ we have
		$$\tr{(\tilde{\A}^{\ell})} = d^{\ell}-a.$$   
	\end{lemma}
	\begin{proof}
		The claim follows easily upon applying Newton's identities to the roots of the irreducible polynomial $2+x+\cdots +x^k$.
	\end{proof}
	
	In particular, since $\tr(\tilde{\A}) \geq 0$, we have $a \leq d$. Hence, an upper bound for the number of orbits is $$w \leq 1 + dk.$$

	Next, we investigate a lower bound for the number of orbits $w$ that will allow us to refine the relationships between $d,k$ and $a$. 
	\begin{lemma}\label{lemm:tr}
		$\tr(\tilde{\A}+\cdots +\tilde{\A}^{ \lfloor k/2 \rfloor}) \leq \lfloor k/2\rfloor w.$
	\end{lemma}
	\begin{proof}
		Set $S = \tr(\tilde{\A}+\cdots +\tilde{\A}^{ \lfloor k/2 \rfloor})$. The trace $\tr(\tilde{\A}^j)$ counts the number of walks of length $j$ starting and ending at the same orbit $\mathcal{O}\in V(\tilde{G})$. Let $j \leq \lfloor k/2 \rfloor $ and suppose that $\pi: \mathcal{O} \to \mathcal{O}_1 \to \ldots  \to \mathcal{O}_j = \mathcal{O}$ and $\pi': \mathcal{O} \to \mathcal{O}_1' \to \ldots \to \mathcal{O}_{j}'=\mathcal{O}$ are two distinct closed walks of length $j$ starting and ending at $\mathcal{O}$. These two walks give rise to two distinct walks of the same length from $v(\mathcal{O})$ to $\mathcal{O}$: 
		$$v(\mathcal{O}) \to v_1 \to \ldots \to v_j=o^a(v(\mathcal{O}))$$
		$$v(\mathcal{O}) \to v_1' \to \ldots \to v_j'=o^b(v(\mathcal{O}))$$
		Then, combining these two walks, we can construct the following two new walks that have length $2j \leq k$, thus reaching a contradiction (see Figure~\ref{fig:1} to visualise this property): 
		$$v(\mathcal{O}) \to v_1 \to ... \to v_j=o^a(v(\mathcal{O})) \to o^a(v_1') \to \ldots \to o^a(o^b(v(\mathcal{O}))) = o^{a+b}(v(\mathcal{O}))$$
		$$v(\mathcal{O}) \to v_1' \to \ldots \to v_j'=o^b(v(\mathcal{O})) \to o^b(v_1) \to \ldots \to o^b(o^a(v(\mathcal{O}))) = o^{a+b}(v(\mathcal{O})).$$
		Now, let $X$ be the set of triplets $(\mathcal{O},j,z)$ with:
		
		\begin{enumerate}
			\item $\mathcal{O}$ an orbit 
			\item $j$ a positive integer with $1 \leq j \leq \lfloor \frac{k}{2}\rfloor$ 
			\item $z$ a node from $\mathcal{O}$ with $d_G(v(\mathcal{O}),z) = j$. 
		\end{enumerate}
		
		On the one hand, for each orbit $\mathcal{O}\in V(\tilde{G})$ we have seen that there cannot exist two distinct walks of length $j$ from $v(\mathcal{O})$ to any other node from $\mathcal{O}$ if $1 \leq j \leq \lfloor k/2\rfloor$. Therefore 
		$$|X| \leq w \lfloor k/2 \rfloor .$$ 
		
		On the other hand, every closed walk in $\tilde{G}$ of length $j$ starting and ending in $\mathcal{O} \in V(\tilde{G})$ corresponds to a walk in $G$ from $v(\mathcal{O})$ to some other node $z \in \mathcal{O}$ of length $j$. Hence $S \leq |X|$. The claim follows by combining these two inequalities. 
		\begin{figure}
			\begin{center}
					\includegraphics[scale=0.75]{figure1-excess.png}
			\caption{The existence of two distinct walks of lengths $\leq k/2$ from $v(r)$ to $o^a(v(r))$ and $o^b(v(r))$ with $a \neq b$ would imply the existence of two distinct walks of length $\leq k$ from $v(r)$ to $o^{a+b}(v(r))$.}
			\label{fig:1}
		\end{center}
	\end{figure}
\end{proof}

We now finally deduce that $(d,k;+1)$-digraphs can exist only for small values of $d$ and $k$. Recall that we can assume that $d \geq 3$ and $k \geq 5$ by the results of~\cite{MilMirSil,Sil,Tuite2024}.

\begin{theorem}\label{thm:main theorem}
	If $G$ is a $(d,k;+1)$-digraph, then $d$ and $k$ satisfy one of the following:
	\begin{itemize}
		\item $k = 5$ and $3 \leq d \leq 11$,
		\item $k = 6$ and $3 \leq d \leq 4$,
		\item $k = 7$ and $3 \leq d \leq 4$, or
		\item $k = 9$ and $d = 3$.
	\end{itemize}
\end{theorem}
\begin{proof} By summing the identity in Lemma~\ref{lem:trace A l} for $1 \leq \ell \leq k-1$ and substituting into Lemma~\ref{lemm:tr}, we deduce that 
	\begin{equation}\label{ineq:1}
		d+\dots +d^{\lfloor k/2 \rfloor}-\lfloor k/2 \rfloor a \leq \lfloor k/2 \rfloor (1+ak). \end{equation}
	As $a \leq d$, this implies that 
	\begin{equation}
		\label{ineq:2}
		d+\dots +d^{\lfloor k/2 \rfloor} \leq \lfloor k/2 \rfloor (1+d+dk). \end{equation}
	It can be easily checked that Inequality (\ref{ineq:1}) can be satisfied only if $d,k$ are one of the listed pairs.
\end{proof}
We will refer to the 14 pairs in Theorem~\ref{thm:main theorem} as \emph{exceptional pairs}. In the next section we rule out the existence of a $(d,k;+1)$-digraph for these exceptional pairs.

\section{Non-existence for small $d$ and $k$}\label{sec:nonexistence}

In~\cite{Tuite2024} a $(d,k;+1)$-digraph is called \emph{minimal} if there is no $(d',k;+1)$-digraph for $2 \leq d' < d$. To prove non-existence of $(d,k;+1)$-digraphs for the exceptional pairs, it suffices to prove that there are no minimal $(d,k;+1)$-digraphs with these parameters. It is known from the same paper that the fixed points of any non-identity automorphism of a minimal $(d,k;+1)$-digraph form one of the following:
\begin{itemize}
	\item the empty set,
	\item a pair of vertices $\{ u,o(u)\} $ that constitutes a transposition in $o$, or
	\item a directed $(k+2)$-cycle.
\end{itemize}
Applying this result to the outlier automorphism itself shows that a minimal $(d,k;+1)$-digraph must have one of the following three structures:
\begin{itemize}
	\item all outlier-orbits of $G$ have the same length (then $G$ is called \emph{outlier-regular}),
	\item the vertices with smallest order $\omega (G)$ induce a directed $(k+2)$-cycle (\emph{Type A} digraphs), or
	\item the outlier automorphism contains a unique transposition, i.e.\ $m_2 = 1$ in the permutation cycle structure of $G$ (\emph{Type B} digraphs).    
\end{itemize}

First we deal with the outlier-regular case. 

\begin{theorem}
	There are no outlier-regular $(d,k;+1)$-digraphs for $d,k \geq 2$.
\end{theorem}
\begin{proof}
	We have $w = \frac{N}{\omega (G)}$ and $a = \frac{w-1}{k}$, both of which must be integers. From Section~\ref{sec:pseudodigraph} we must also have \[ d+\dots +d^{\lfloor k/2 \rfloor}-\lfloor k/2 \rfloor a \leq \lfloor k/2 \rfloor (1+ak) \text{ and } a \leq d,\]
	and it is straightforward to check that none of the exceptional pairs satisfy these conditions.
\end{proof}

We now proceed to rule out the Type A and B cases. An important fact is that for such digraphs there can be at most three distinct outlier orders.

\begin{lemma}\label{lem:at most 3 cycle lengths}
	Any minimal $(d,k;+1)$-digraph $G$ has at most three distinct vertex orders, and if $G$ has vertex orders $s_1 < s_2 < s_3$, then $s_1 = 2$ and $s_2 = k+2$.
\end{lemma}
\begin{proof}
	Suppose that the outlier automorphism contains cycles of lengths $\ell _1 < \ell _2 < \dots < \ell r$. Suppose that $\ell _i > 2$, $1 \leq i \leq r$ and consider the automorphism $o^{\ell _i}$. This fixes at least three vertices, but is not equal to the identity, and so its fixset must be a $(k+2)$-cycle $C$ in $G$. If $v$ is a vertex on $C$, then $o(v)$ must be the vertex preceding $v$ on $C$. Therefore these fixed points correspond to a cycle of length $\ell _i = k+2$ in $o$. As this applies to any non-transposition cycle in $o$, the result follows. 
\end{proof}
In particular, Lemma~\ref{lem:at most 3 cycle lengths} shows that for Type A minimal $(d,k;+1)$-digraphs there are exactly two distinct cycle lengths, the smaller of which is $\omega (G) = k+2$. 

\begin{theorem}
	There are no minimal Type A or Type B digraphs for the exceptional pairs.
\end{theorem}
\begin{proof}
	Let $G$ be a minimal $(d,k;+1)$-digraph with $(d,k)$ belonging to one of the exceptional pairs from Theorem~\ref{thm:main theorem}.
	
	Suppose first that $G$ is Type A with a $(k+2)$-cycle $C$. Let $u \rightarrow u'$ be an arc of $C$. The automorphism $o^{k+2}$ fixes $C$, but has no other fixed points; in particular, it must permute the vertices of $N^+(u)-\{ u'\} $ without any fixed points. Suppose that the shortest cycle of the permutation $o^{k+2}$ on $N^+(u)-\{ u'\} $ has length $t$; then the automorphism $o^{(k+2)t}$ must fix every vertex of $G$, and it follows that every cycle of $o^{k+2}$ on $N^+(u)-\{ u'\} $ has length $t$, and so $t \mid d-1$. It now follows from Lemma~\ref{lem:at most 3 cycle lengths} that the cycle length of $o$ are $k+2$ and some divisor $\ell > k+2$ of $(k+2)(d-1)$. Therefore \[ 1+m_{\ell } = w, a = \frac{N-k-2}{k\ell }. \] We only need to examine cases in which $a \leq d$, so we can assume that \[ N \leq k+2+k\ell d \leq k+2+k(k+2)(d-1)d.\] This inequality does not hold for any of the exceptional pairs.
	
	Now suppose that $G$ is a Type B digraph. Let $u,o(u)$ be the two vertices in the unique transposition in $o$. Consider the action of $o^2$ on the vertices of $N^+(u) = \{ u_1,\dots ,u_d\} $. As $o^2$ fixes only the vertices $u,o(u)$, $o^2$ does not fix any vertices of $N^+(u)$. Consider the cycle structure of the action of $o^2$ on $\{ u_1,\dots ,u_d\} $. If this action has cycles of lengths $s,t$, where $s < t$, then the automorphism $o^{2s}$ must fix a $(k+2)$-cycle, which is impossible. Therefore, $o^2$ acts on $\{ u_1,\dots ,u_d\} $ as a collection of $t$-cycles for some $2 \leq t \leq d$. Then the automorphism $o^{2t}$ fixes all vertices of $\{ u,o(u)\} \cup N^+(u)$, and so all vertex orders must divide $2t$.
	
	If $G$ has three distinct cycle lengths $2 < \ell _2 < \ell _3$, where $\ell _2$ and $\ell _3$ both divide $2t$, then by Lemma~\ref{lem:at most 3 cycle lengths} we have $\ell _2 = k+2$. It is easily checked that this is impossible for the exceptional pairs, and so in each of these cases each vertex of $V(G) -\{ u,o(u)\} $ has order $\ell _2$, where $\ell _2 > 3$, $\ell _2 \mid 2t$, $t \mid d$ and $t > 1$. We have $1+m_{\ell } = w$, so $a = \frac{m_{\ell }}{k} = \frac{N-2}{\ell k}$. Therefore, as we need only check cases when $a \leq d$, we can assume that $N \leq 2+d\ell k \leq 2+d^2k$, which does not hold for the exceptional cases. 
\end{proof}

This concludes the proof of our main result.

\begin{theorem}
	There are no $(d,k;+1)$-digraphs for $d,k \geq 2$.
\end{theorem}

\section{Conclusion}\label{sec:conclusion}

In this article we have proven the non-existence of directed graphs with excess one. We now briefly explore the implications of this result for the study of networks close to the Moore bound. 

In~\cite{Baskoro} it is proved that any almost Moore digraph $G$ without self-repeats always contains a subdigraph $H$ which is itself an almost Moore digraph, whose permutation cycle structure is constrained in one of the following forms:

\begin{itemize}
	\item All orbits with respect to the repeat function $r$ have the same length $s$, so that the order of $H$ is $N = s m_s$, where $m_s$ is the number of orbits of length $s$.
	\item A mixture of orbits with respect to the repeat function $r$ of length $s$ and $2s$, so that the order of $H$ is $N = sm_s + 2sm_{2s}$, with $m_s$ orbits of length $s$ and $m_{2s}$ orbits of length $2s$.
\end{itemize}

However, in the second case, the subdigraph induced by the nodes belonging to cycles of length $s$ is a digraph of excess one. Since we have already established that digraphs of excess one do not exist, this case cannot occur. 

Consequently, to prove the non-existence of almost Moore digraphs in general, it is sufficient to focus only on those almost Moore digraphs whose subdigraphs have all orbits of the same length. In other words, the non-existence of digraphs of excess one allows us to reduce the general problem to analysing almost Moore digraphs with a very restricted and simpler cycle structure.

The next natural question is to classify the pairs $(d,k)$ for which there exists a $(d,k;+2)$-digraph. It is known that there are two such digraphs for $(d,k) = (2,2)$, but none with $d = 2$ and $k \geq 3$, see~\cite{d2e2nondiregular,d2e2diregular}. We make the following conjecture.

\begin{conjecture}
	For $d,k \geq 2$ there exists a $(d,k;+2)$-digraph only for $(d,k) = (2,2)$.
\end{conjecture}

It is also interesting to note that, whilst there are no non-trivial digraphs with excess one, there do exist non-trivial almost Moore digraphs. This suggests that it is easier to approach the Moore bound from below than from above, and inspires the following question.

\begin{problem}\label{problem}
	Is there a pair $(d,k)$ for which the smallest possible defect of a $d$-out-regular digraph with diameter $k$ is less than the smallest possible excess of a $d$-out-regular $k$-geodetic digraph? If yes, what can be said about the proportion of such pairs?   
\end{problem}
We note that there do exist such pairs for the undirected version of Problem~\ref{problem}: the largest graph with degree $6$ and diameter $2$ has defect $5$, whereas the smallest graph with degree $6$ and girth $5$ has excess $3$.

\section*{Acknowledgments}
This research was carried out during visits by James Tuite to the University of Primorska and Universitat de Lleida, and he thanks both institutions for the generous hospitality. Funding for the former visit was partially provided by a grant from the OU Crowther Fund. The authors thank Grahame Erskine for pointing out the example for undirected graphs after Problem~\ref{problem}. The second and third authors also gratefully acknowledge the support provided by the consolidated research group 2021SGR-00434.

\end{document}